\documentclass[12pt]{article}
\usepackage{graphicx} 
\usepackage{amsmath,amssymb,enumerate,xcolor}
\usepackage{float}

\title{On the  existence of   geodesic vector fields   on closed surfaces}
\author{   Vladimir S.\ Matveev\footnote{
Institut f\"ur Mathematik, Friedrich Schiller Universit\"at Jena,
07737 Jena,  Germany  \ \ \quad {\tt  vladimir.matveev@uni-jena.de}}}

 \date{ }

\newtheorem{defin}{Definition}

\newtheorem{example}{Example}

\newcommand{\weg}[1]{}

\usepackage{biblatex} 
\begin{document}

\maketitle
\begin{abstract} 
We construct an example of a Riemannian   metric   on the 2-torus such   that its universal cover does not admit global  Riemann normal  coordinates.

\vspace{1ex}

\vspace{1ex}
\end{abstract}

\section{Introduction}

\begin{defin}   We call a vector field $v =v(x)$  on a Riemannian manifold $(M^n, g)$  {\rm geodesic}, if its length is identically  $1$ and if $\nabla^g_vv=0, $  where  $\nabla^g$ is the Levi-Civita connection of $g$.
\end{defin}

 Clearly, a vector field is geodesic if and only if  any orbit of its flow is an arc-length parameterised geodesic.
 
\begin{example} For the metric 
\begin{equation}\label{eq:1}
g= (dx^1)^2 + \sum_{i,j=2}^n  g_{ij}(x) dx^idx^j \, ,\end{equation}
 the  vector field 
$\frac{\partial}{\partial x^1} $ is geodesic. 
\end{example}

In dimension two the formula  \eqref{eq:1} reads  
\begin{equation}\label{eq:2}
g= dx^2 +   f(x,y) dy^2. \end{equation}
Coordinates such that the metric has the form \eqref{eq:2}
 are called {\it Riemann normal} coordinates. It is known that, in dimension two,  for any geodesic vector field  there exists a local coordinate system $(x,y)$  such that the metric has the form  \eqref{eq:2} and the  vector field is $\frac{\partial}{\partial x}$.  
 
The goal of this paper is to construct an example of a Riemannian two-torus  $(T^2, g)$ such that its universal cover $(\mathbb{R}^2, \tilde g),  $ where $\tilde g$ denotes the lift of $g$,   has no geodesic vector field.      Any sufficiently  small $C^2$-perturbation  of   this metric has the same property.  The example can be easily generalised to closed surfaces of negative Euler characteristic.

We have the following two  motivations for studying the problem. The first one is related to the very recent paper \cite{MP} studying conformal product structures on K\"ahler manifolds.  \cite[Corollary 4.6]{MP} guarantees 
the existence of a  geodesic vector field on compact K\"ahler manifold of real dimension $n \ge  4$ carrying
an orientable conformal product structure with non-identically zero Lee form.    \cite[Proposition 4.7]{MP} uses the results of the present paper to show the existence of direct product compact  K\"ahler metrics   with  no orientable conformal
product structure with non-identically zero Lee form.

Another motivation comes from the theory of integrable geodesic flows on closed surfaces.  \cite[Theorem 1.6]{Bialy10} implies that  for any Riemannian 2-torus $(T^2, g)$ such that the geodesic flow is integrable and the integral satisfies $\aleph$-condition, 
see \cite[Definition 1.3]{Bialy10}, there exists a geodesic vector field on the universal cover $(\mathbb{R}^2, \tilde g)$, where $\tilde g$ denotes the lift of $g$.  Our example is an ``easy to construct'' example of $\aleph$-nonintegrable geodesic flow.   Recall that though generic geodesic flow is not integrable,  proving that a geodesic flow is nonintegrable or constructing  an example of an nonintegrable geodesic flow is not an easy task, see e.g. \cite[\S 10]{Burns} and \cite[\S 3]{BMMT}.

\subsection*{Acknowledgement. } I  warmly thank Andrei  Moroianu for asking the question that lead to this paper, and for encouraging me to write the paper. 
I was supported by the DFG (projects 455806247 and 529233771)  and  by the ARC Discovery Programme DP210100951. 

\section{ Example and proof of nonexistence of geodesic vector field }

\begin{figure}[H]
\begin{minipage}{.45\textwidth} Take the standard sphere with the standard metric. Next,  take a small $\varepsilon>0 $ and change the topology of the manifold in the $\varepsilon$-neighborhood of the south  pole      by gluing  
 a handle  in the neighborhood. The metric outside  the neighborhood is not changed, the metric in the modified neighborhood can be chosen arbitrary such that the obtained metric on the two-torus is smooth, see Fig. \ref{fig:1}.

 \end{minipage}\begin{minipage}{.05\textwidth} {\color{white} .}  \end{minipage}\begin{minipage}{.4\textwidth}  \centering
     \includegraphics[width=0.37\linewidth]{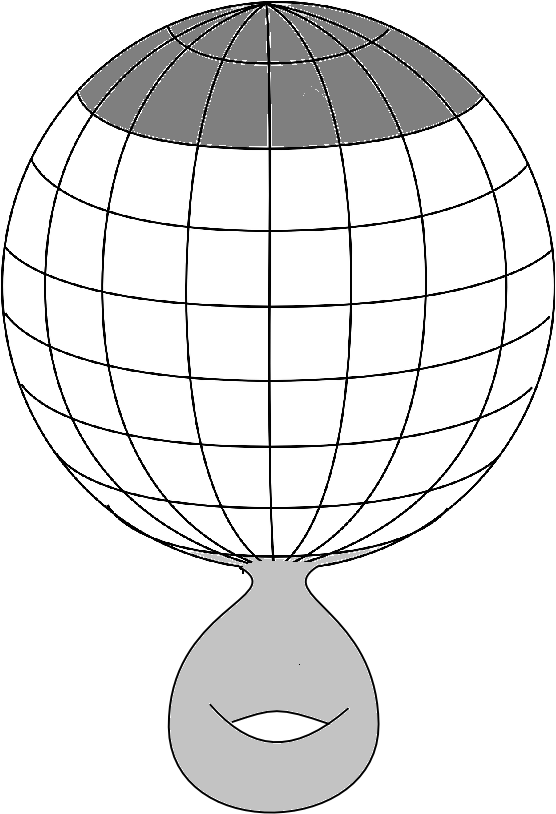}
     \caption{The torus made of the   sphere: the dark-gray part is the $2\varepsilon$-ball around the north pole. The surgery was made in the light-gray part. }
     \label{fig:1}\end{minipage}
 \end{figure}
  
 We consider the universal cover $\mathbb{R}^2$ and denote by $\tilde g$ the lift of the metric. 
Let us show  that  $(\mathbb{R}^2, \tilde g)$ does not admit a geodesic vector field. We assume it does, denote the geodesic vector field by $v$, and find a contradiction.

In order to do it, consider the  circle of radius  $2\varepsilon$ around  the north pole of the initial sphere and consider one of its lifts 
 $C_{2\varepsilon}(N_0)= \partial B_{2\varepsilon}(N_0)$.  
Let us show that 
our geodesic  vector field $v$ is necessary transversal to it. Arguing by  the method of contradiction, assume  there exists a point where the vector field $v$  is tangent to the  circle. Consider  the geodesic $\gamma$ starting from this point in the direction of $v$. This geodesic, and also geodesics close to $\gamma$,  
do not enter the ``light gray'' region where we changed the sphere. Therefore, any geodesic  $ \gamma_1$ starting from a nearby point in the direction of our vector field  intersects $\gamma$, as any two geodesics on the sphere intersect each other. This    gives a contradiction, since  velocity vectors of both geodesics at the point of intersection should be $v$.

Thus, our distribution is transversal   to the circle at every point.  Then, the index  of the restriction of   $v$   to $B_{2 \varepsilon}(N_0)$ is nonzero.  But the index  must be zero since $v$ is never zero. 
The contradiction proves the nonexistence of a geodesic vector field. 

Note also that  in the proof we used the following properties  of the standard metric of the sphere only:
\begin{enumerate}

\item The geodesic starting at a point of  the $2\varepsilon$- circle and tangent to it does not reach the $\varepsilon$-neighborhood of the south pole within the time $2\pi$.

\item  Two geodesics    always intersect. 
 
\end{enumerate}
 
 These properties  are evidently fulfilled for any sufficiently small perturbation, in the $C^2$-topology, of the standard metric of the sphere.  This implies that  one can construct such an example in the  real-analytic category. Moreover, by attaching  more than one handle  in the ``light gray'' region  one can construct an example of a closed Riemannian surface   of arbitrary negative Euler characteristic such that the universal cover does not admit a geodesic vector  field.

\printbibliography
\end{document}